\documentclass[11pt]{amsart}
\usepackage{amsmath, amsthm, latexsym, amssymb, amsfonts, epsfig, epsf}
\usepackage[dvips]{color}

\newenvironment{pf}{\proof[\proofname]}{\endproof}
\theoremstyle{plain}
\newtheorem{Th}{Theorem}[section]
\newtheorem{Cor}[Th]{Corollary}
\newtheorem{Prop}[Th]{Proposition}

\numberwithin{equation}{section} \theoremstyle{definition}
\newtheorem{Rem}[Th]{Remark}

\newtheorem{Ex}{Example}

\newcommand{\cal}[1]{\mathcal{#1}}
\newcommand{\F}{\mathbb F}
\newcommand{\T}{(\mathbb F_q^*)^n}

\newcommand{\N}{\mathbb N}
\newcommand{\Z}{\mathbb Z}
\newcommand{\R}{\mathbb R}
\newcommand{\Pp}{\mathbb P}

\newcommand{\cC}{\cal C}

\newcommand{\cL}{\cal L}
\newcommand{\cO}{\cal O}

\newcommand{\cP}{\cal P}

\newcommand{\spn}{\operatorname{span}}

\newcommand{\chull}{\operatorname{convex\,hull}}

\newcommand{\rs}[1]{Section~\ref{S:#1}}

\newcommand{\rp}[1]{Proposition~\ref{P:#1}}

\newcommand{\rr}[1]{Remark~\ref{R:#1}}
\newcommand{\rex}[1]{Example~\ref{ex:#1}}
\newcommand{\re}[1]{(\ref{e:#1})}
\newcommand{\rc}[1]{Corollary~\ref{C:#1}}
\newcommand{\rt}[1] {Theorem~\ref{T:#1}}

\newcommand{\rf}[1]{Figure~\ref{F:#1}}

\begin{document}

\author{Ivan Soprunov, Jenya Soprunova}
\title{Bringing Toric Codes to the next dimension}
\maketitle

\begin{abstract} This paper is concerned with the minimum distance
computation for higher dimensional toric codes defined by
lattice polytopes in $\R^n$. We show that the minimum distance is multiplicative 
with respect to taking the product of polytopes, and  behaves in a simple way
when one builds a $k$-dilate of a pyramid over a polytope. This allows us to 
construct a large class of examples of higher dimensional toric codes 
where we can compute the minimum distance explicitly. 
\end{abstract}


\section*{Introduction}

In their fundamental work ``Algebraic Geometric codes'' \cite{TV}, 
Tsfasman and Vl\u adu\c t proposed a general framework for applying higher-dimensional 
algebraic varieties to constructing error-correction codes.
Toric varieties provide a
rich class of examples of algebraic varieties where one can make explicit calculations 
using the connection of toric varieties to geometry of lattice polytopes.

In general, given a normal projective variety $X$ over a finite field $\F_q$,
a line bundle~$L$ on~$X$  defined over $\F_q$, and a finite set $S=\{p_1,\dots,p_s\}$ 
of $\F_q$-rational points on~$X$, one defines an {\it evaluation code} $\cC(X,L,S)$
by evaluating global sections of $L$ at the points of~$S$. In other words, it is the image of the map
$$ev:\Gamma(X,L)\to\F_q^s,\quad f\mapsto(f(p_1),\dots,f(p_s)).$$

Now let $P$ be a lattice polytope in $\R^n$. It defines a projective toric variety~$X_P$ 
and an ample line bundle $L_P$ (see \cite{F}, Section 3.4). If for $S$ we take the algebraic
torus $\mathbb T=(\F_q^*)^n$, the above construction gives a {\it toric code} 
$\cC_P=\cC(X_P,L_P,\mathbb T)$. These codes were first
introduced by Hansen for $n=2$ (see \cite{Ha1, Ha2}) and subsequently studied by Joyner, Little,
Schenck, and others \cite{Jo, LiSche, LiSchw, Ru, SoSo}.
Since the global sections $\Gamma(X_P,L_P)$ can be identified with polynomials
whose monomials are the lattice points in $P$ (see \cite{F}, Section 3.4), 
a toric code is completely determined by~$P$ (see the definition we give in \rs{toric}).

Therefore, one should expect parameters of a toric code, such as the dimension and
the minimum distance, to be expressible in terms of
the polytope~$P$. Indeed, when $ev$ is injective (this is true, for example, when $P$
lies in the $n$-cube $[0,q-2]^n$) the dimension of the code is the number of lattice points 
$|P\cap\Z^n|$ in~$P$. More about injectivity can be found in \cite{Ru}. 

As for the minimum distance of $\cC_P$, the following results reveal its strong connection with geometry of~$P$.  First, for $n=2$ and large enough~$q$, we were able to write bounds for 
the minimum distance in terms of the {\it Minkowski length} of $P$ (see~\cite{SoSo}). 
Second, there are explicit formulas for the minimum distance of the toric codes $\cC_{k\,\Delta}$
and  $\cC_{\Pi}$, where the polytopes are the $k$-dilate of the standard $n$-simplex $\Delta$ and
a rectangular box $\Pi=[0,a_1]\times\cdots\times [0,a_n]$, respectively.
 (See~\cite{LiSchw} for an elementary treatment of these cases, 
 also see Corollaries \ref{C:simplex} and \ref{C:box}.)
 The corresponding toric varieties
and line bundles are $\Pp^n$ and $\cO(k)$ in the first case, and
$\Pp^1\times\cdots\times\Pp^1$ and $\cO(a_1)\times\dots\times\cO(a_n)$ in the second case.

The goal of the present paper is to construct a larger class of examples where we can compute
the minimum distance explicitly. We show that the minimum distance is multiplicative 
with respect to taking the product of polytopes, and  behaves in a simple way
when one builds a $k$-dilate of a pyramid over a polytope (\rt{product}, \rt{pyramid}).
Thus, starting with a known toric code we can compute the minimum distance for
a toric code in higher dimensions. For example, let $P$ be obtained from a lattice segment by
a sequence of operations, each of them being either
a multiplication by a segment or taking a pyramid and dilating. Then we get
an explicit formula for the minimum distance of the corresponding toric code $\cC_P$ (\rt{mix}).

The methods we use are of combinatorial nature and our approach is influenced by
Serre's way of finding the maximal number of $\F_q$-zeroes of homogeneous polynomials
of fixed total degree (see \cite{Serre}). A similar idea also appears in Hansen's work on 
toric surfaces \cite{Ha1,Ha2} and in \cite{Ru}, but the techniques used there rely heavily on 
intersection theory on toric surfaces. In \rs{param} we consider the asymptotics of the parameters of toric
 codes constructed in \rt{mix}. We conclude the paper with a few examples of 2- and  3-dimensional toric codes.

We thank the anonymous referees for their comments and suggestions that helped
improve the paper.
 

\section{Toric codes}\label{S:toric}

Let $q$ be a prime power and  $P$ be a lattice polytope in
$\R^n$ contained in the $n$-cube $K_{q}^n=[0,q-2]^n$. We will
denote by $\cL(P)$ the vector space over the finite field~$\F_q$
spanned by the monomials $x^m=x_1^{m_1}\cdots x_n^{m_n}$ corresponding to the lattice points 
$m=(m_1,\dots,m_n)$ of $P$,
$$\cL(P)=\spn\{\,x^m\ |\ m\in P\cap\Z^n\}.$$

The {\it toric code} $\cC_P$ is a linear code of block length $N=(q-1)^n$
whose codewords are the vectors of values
of $f\in\cL(P)$ at all points of the algebraic torus $(\F_q^{\,*})^n$ (in some linear order):
$$\cC_P=\{\left(f(\xi), \xi\in(\F_q^{\,*})^n\right)\ |\ f\in\cL(P)\}.$$
Note that the vectors obtained by evaluating the monomials $\{\,x^m\ |\ m\in P\cap\Z^n\}$
are linearly independent over $\F_q$, since $P\subseteq K_{q}^n$ (see \cite{Ru}). Therefore, 
$\cC_P$ has dimension $|P\,\cap\,\Z^n|$, the number of lattice points in $P$.

 The
{\it weight} of each nonzero codeword equals the number of points
$\xi\in(\F_q^{\,*})^n$ where the corresponding polynomial does not vanish.
We will denote it by $w(f)$.

Let $Z(f)$ denote the number of points in $(\F_q^{\,*})^n$ where $f$ vanishes,
these points will be called $\F_q$-zeroes or just zeroes of $f$.
Then the minimum distance $d(\cC_P)$, which is also the minimum
weight, equals
$$d(\cC_P)=(q-1)^n-\max_{0\neq f\in\cL(P)}Z(f).$$

Throughout the paper we let $Z_P$ denote the maximum number of zeroes 
over all nonzero $f\in\cL(P)$. Obviously, $Z_P=(q-1)^n-d(\cC_P)$.



\section{Main Theorems}

\subsection{Product of polytopes}
Let $P$ and $Q$ be lattice polytopes in $K_q^m$ and $K_q^n$, respectively. 
Then the product $P\times Q$ is a lattice polytope in 
$K_{q}^{m+n}$.  In the following theorem we prove that the minimum distance for the toric code
$C_{P\times Q}$ is simply the product of the minimum distances of $\cC_P$ and $\cC_Q$.

\begin{Th}\label{T:product} Let $P\subseteq K_q^m$ and $Q\subseteq K_q^n$ be lattice polytopes.
Then  
$$d(\cC_{P\times Q})=d(\cC_P)\,d(\cC_{Q}).$$
\end{Th}

\begin{pf} First, we will show that there exists a polynomial $f\in\cL(P\times Q)$ whose weight
is exactly $d(\cC_P)\,d(\cC_{Q})$. Let $g\in\cL(P)$ be a polynomial with 
$w(g)=d(\cC_{P})$ and let $h\in\cL(Q)$ be a polynomial with 
$w(h)=d(\cC_{Q})$.
Consider the polynomial 
$$f(x_1,\dots,x_{m},y_1,\dots,y_n)=g(x_1,\dots, x_{m})\,h(y_1,\dots, y_{n}).$$
Clearly $f$ lies in  $\cL(P\times Q)$. Note that $f(\xi_1,\dots, \xi_m,\eta_1,\dots,\eta_n)$
is non-zero if and only if both $g(\xi_1,\dots, \xi_{m})$ and $h(\eta_1,\dots, \eta_{n})$ are non-zero.
Therefore, the weight of $f$ equals the product of the weights of $g$ and $h$, i.e.,
 $w(f)=d(\cC_P)\,d(\cC_{Q})$.

 Now consider an arbitrary nonzero  $f\in\cL(P\times Q)$.  We need to show that
  $w(f)\geq d(\cC_P)\,d(\cC_{Q})$. We can write
\begin{equation}\label{e:goodform}
f(x_1,\dots,x_{m},y_1,\dots,y_n)=\sum_{u\in Q\cap\Z^n}f_u(x_1,\dots,x_{m})y^u,
\end{equation}
for some polynomials $f_u$. Note that each $f_u$ is supported in $P$.

Consider an affine subspace 
$L_\xi=\{(\xi,y)\in (\F_q^*)^{m+n}\}$ through a point $\xi\in  (\F_q^*)^{m}$.
Let  $s$ denote the number of those $L_\xi$ where $f$ is identically zero.  
Then the number of zeroes of $f$ in $(\F_q^*)^{m+n}$ is bounded from above by
\begin{equation}\label{e:long}\nonumber
Z(f)\leq (q-1)^ns+Z_Q\left((q-1)^{m}-s\right)
=Z_Q(q-1)^{m}+s\left((q-1)^n-Z_Q\right).
\end{equation}
Indeed, on every $L_\xi$ where $f$ is identically zero,
$f$ has exactly $(q-1)^n$ zeroes. Furthermore, on every $L_{\xi}$ where $f$ is not
identically zero, it has at most $Z_Q$ zeroes since the
polynomial obtained by substituting $x_i=\xi_i$, for $1\leq i\leq {m}$
is nonzero and supported in~$Q$ (see \re{goodform}).

Now we will obtain a bound for $s$. 
The fact that $f$ vanishes identically on  $L_\xi$ implies that
$\xi$ is a common zero of the polynomials $f_u$ from \re{goodform}, hence $s$ is the
number of common zeroes of the $f_u$ in
$(\F_q^*)^{m}$. But the number of common zeroes of the
$f_u$ is at most the number of zeroes of each (nonzero) $f_u$. Therefore $s\leq Z(f_u)\leq Z_P$,
since $f_u$ is supported in $P$. We obtain
\begin{equation}\nonumber
Z(f)\leq Z_Q(q-1)^{m}+Z_P((q-1)^n-Z_Q)
=(q-1)^{m+n}-d(\cC_{P})\,d(\cC_{Q}).
\end{equation}
Therefore, the weight of $f$  is at least $d(\cC_{P})\,d(\cC_{Q})$, which completes the proof. 
 \end{pf}

\subsection{Pyramids over polytopes}

Let $Q$ be an $n$-dimensional lattice polytope contained in $K_{q}^n$.
Furthermore,  let $\{kQ\ |\ 0\leq k\leq N\}$ 
be a sequence of $k$-dilates of $Q$,  contained in $K_{q}^n$,
for some $N<q-1$. It defines the sequence of minimum distances $d(\cC_{kQ})$ for the corresponding
toric codes. It turns out that this sequence decreases sufficiently quickly, as the following proposition shows.

\begin{Prop}\label{P:decrease} 
Let $Q$ be a polytope of $\dim Q\geq 1$, and let $\{kQ\ |\ 0\leq k\leq N\}$ 
be a sequence of $k$-dilates of $Q$,  contained in $K_{q}^n$. Then the 
sequence 
$d_k=d(\cC_{kQ})$ satisfies 
\begin{equation}\label{e:decrease}
\frac{d_k}{d_{k-l}}\leq 1-\frac{l}{q-1},\quad\text{for all }\ 0\leq l\leq k.
\end{equation}
\end{Prop}

\begin{pf} To shorten our notation we let $Z_k=Z_{kQ}=(q-1)^n-d_k$. Then \re{decrease}
is equivalent to 
\begin{equation}\label{e:equiv}
Z_k\geq Z_{k-l}+l(q-1)^{n-1}-\frac{lZ_{k-l}}{q-1}.
\end{equation}
This is what we are going to show.
Consider a polynomial $g\in\cL\left((k-l)Q\right)$ with $Z(g)=Z_{k-l}$. 
Since $Q$ is at least 1-dimensional,
it contains a primitive lattice segment $I\subseteq Q$. Then, for any $0\leq l\leq k$, we have
\begin{equation}\label{e:include}
kQ=(k-l)Q+lQ\supseteq (k-l)Q+lI.
\end{equation}
Since unimodular equivalent polytopes result in monomially equivalent toric codes (see \cite{LiSchw}) 
we may assume that $I$ is the unit segment $I=[0,e_n]$, where $e_n=(0,\dots,0,1)$.
Then \re{include} implies that the polynomial
$$f(x_1,\dots,x_n)=g(x_1,\dots,x_n)\prod_{j=1}^l(x_n-a_j)$$
is contained in $\cL(kQ)$ for any choice of $a_1,\dots,a_l\in\F_q^*$.
We have
\begin{equation}\label{e:almost}
Z_k\geq Z(f)=Z(g)+l(q-1)^{n-1}-\sum_{j=1}^lZ(g|_{x_n=a_j}).
\end{equation}

Notice that $\T$ is the union of hyperplanes $\{x_n=a\}$, for $a\in\F_q^*$. Therefore,
$\sum_{a\in\F_q^*}Z(g|_{x_n=a})=Z(g)$. We choose $a_1,\dots,a_l\in\F_q^*$ such that
$$\{Z(g|_{x_n=a_j})\ |\ j=1,\dots,l\}$$ are the $l$ smallest integers among the $q-1$ integers
$\{Z(g|_{x_n=a})\ |\ a\in\F_q^*\}$. Since their average is no greater
than the average of the $q-1$ integers, we have:
$$\frac{1}{l}\sum_{j=1}^lZ(g|_{x_n=a_j})\leq \frac{1}{q-1}\sum_{a\in\F_q^*}Z(g|_{x_n=a})=\frac{Z(g)}{q-1}.$$
Combining this inequality with \re{almost} and using $Z(g)=Z_{k-l}$ we obtain the desired
inequality \re{equiv}.

\end{pf}

\begin{Rem}\label{R:stronger}
It follows from the above proof  that if $Q$ contains a lattice segment
of (lattice) length $\lambda$ then the sequence $d_k$ possesses a stronger decreasing property:
$$\frac{d_k}{d_{k-l}}\leq 1-\frac{\lambda l}{q-1},\quad\text{for all }\ 0\leq l\leq k.$$
\end{Rem}

Let $Q$ be an $n$-dimensional lattice polytope contained in $K_{q}^n$.
We denote by $\cP(Q)$ the unit pyramid over $Q$, i.e. 
$$\cP(Q)=\chull\{ e_{n+1}, (x,0)\ |\ x\in Q\}\subseteq K_q^{n+1}.$$
As before a sequence  $\{kQ\ |\ 0\leq k\leq N\}$ 
of $k$-dilates of $Q$,  contained in $K_{q}^n$,
defines the sequence of minimum distances $d(\cC_{kQ})$ for the corresponding
toric codes. In the next theorem we relate 
 $d(\cC_{kQ})$ and $d(\cC_{k\cP(Q)})$.

\begin{Th}\label{T:pyramid}
Let $Q$ be a lattice polytope of $\dim Q\geq 1$, and let $\{kQ\ |\ 0\leq k\leq N\}$ 
be a sequence of $k$-dilates of $Q$,  contained in $K_{q}^n$. Then 
$$d(\cC_{k\cP(Q)})=(q-1)\,d(\cC_{kQ}).$$
\end{Th}

\begin{pf} The proof is similar to the proof of \rt{product}. First, let $g\in\cL(kQ)$ be a polynomial
with weight $d(\cC_{kQ})$ and number of zeroes $Z(g)=(q-1)^{n}-d(\cC_{kQ})$. 
Then, as an element of $\cL(k\cP(Q))$, the polynomial $g$
has exactly $(q-1)Z(g)$ zeroes in $(\F_q^*)^{n+1}$, hence, has weight $(q-1)\,d(\cC_{kQ})$.

 Now we consider an arbitrary nonzero  $f\in\cL(k\cP(Q))$ and show it has weight
 at least $(q-1)\,d(\cC_{kQ})$. We write
\begin{equation}\label{e:goodform2}
f(x_1,\dots,x_{n+1})=\sum_{i=0}^{k}f_i(x_1,\dots,x_{n})x_{n+1}^i.
\end{equation}
for some polynomials $f_i$  supported in (a translation of)
$(k-i)Q$. Let $l$ be the largest value of $i$ for
which $f_i$ is not identically zero. Then \re{goodform2} becomes
\begin{equation}\nonumber
f(x_1,\dots,x_{n+1})=\sum_{i=0}^{l}f_i(x_1,\dots,x_{n})x_{n+1}^i.
\end{equation}

Similar to the proof of \rt{product}, we let $L_\xi$ be the vertical line through a point $\xi\in \T$.
On every $L_\xi$ where $f$ is identically zero,
$f$ has exactly $q-1$ zeroes, and on every $L_{\xi}$ where $f$ is not
identically zero, it has at most $l$ zeroes, since the univariate polynomial
obtained by substituting $x_i=\xi_i$, for $1\leq i\leq {n}$, 
has degree at most $l$.
Then the number of zeroes of $f$ in $(\F_q^*)^{n+1}$  is bounded by
\begin{equation}\label{e:long2}
Z(f)\leq (q-1)s+l\left((q-1)^{n}-s\right)
=l(q-1)^{n}+s\left(q-1-l\right),
\end{equation}
where $s$ is the number of the lines $L_{\xi}$ where $f$ is identically zero.
Recall that $s$ is, in fact, the number of common zeroes of the $f_i$ in
$(\F_q^*)^{n}$, and is at most the number of zeroes of each $f_i$. 
In particular, $s\leq Z(f_l)$. Since $f_l$ is supported in a (vertical) translation
of $(k-l)Q$, we have 
$$s\leq Z(f_l)\leq (q-1)^n-d(\cC_{(k-l)Q})=(q-1)^n-d_{k-l},$$ 
where we abbreviate  $d_k=d(\cC_{kQ})$.
Therefore, \re{long2} implies
\begin{equation}\nonumber
Z(f)\leq l(q-1)^{n}+((q-1)^n-d_{k-l})(q-1-l)
=(q-1)^{n+1}-d_{k-l}(q-1-l).
\end{equation}
Finally, by \rp{decrease}, $d_{k-l}(q-1-l)\geq d_{k}(q-1)$, and hence,
$$Z(f)\leq (q-1)^{n+1}-d_{k}(q-1),$$
which means that the weight of $f$  is at least $(q-1)\,d_k=(q-1)\,d(\cC_{kQ})$.
\end{pf} 

Applying \rt{pyramid} and \rp{decrease} with $l=k-1$ we obtain an upper bound for
the minimum distance of $\cC_{k\cP(Q)}$ in terms of the minimum distance of $\cC_{Q}$:

\begin{Cor}
Let $Q$ be a lattice polytope of $\dim Q\geq 1$, and let $k\cP(Q)$ be a $k$-dilate
of the unit pyramid over $Q$, contained in $K_{q}^{n+1}$. Then 
$$d(\cC_{k\cP(Q)})\leq(q-k)\,d(\cC_{Q}).$$
\end{Cor}

\begin{Rem}\label{R:doublepyr} Using essentially the same arguments as above
one can prove a similar result  for the {\it double pyramid} over $Q$.
Let $\cal{DP}(Q)$ denote the double pyramid over~$Q$, translated so that it is contained in $K_q^{n+1}$:
$$\cal{DP}(Q)=\chull\{e_{n+1},-e_{n+1}, (x,0)\ |\ x\in Q\}+(0,\dots,0,1).$$
Then if $Q$ contains a lattice segment of lattice length 2 then 
$$d(\cC_{k\cal{DP}(Q)})=(q-1)\,d(\cC_{kQ}).$$
Indeed, the proof is exactly the same as in the above theorem, except $l$ should be replaced with $2l$
and \rr{stronger} should be used in place of \rp{decrease}.
\end{Rem}

We conclude this section with a remark about the dimension of the toric code. As we mentioned in \rs{toric}, the dimension of $\cC_P$ equals the number of
lattice points in $P$. It is well known how the number of lattice points behaves when one
takes the product of polytopes or builds a $k$-dilate of a pyramid over a polytope.
Clearly $|(P\times Q)\cap Z^{m+n}|=|P\cap\Z^m||Q\cap\Z^n|$ for lattice polytopes $P\subset \R^m$,
$Q\subset\R^n$. Hence, just as for the minimum distance 
$$\dim(\cC_{P\times Q})=\dim(\cC_P)\dim(\cC_Q).$$
Now if $\cP(Q)\subset \R^{n+1}$ is the pyramid over $Q\subset\R^n$. 
Then the set of lattice points in $k\hspace{1pt}\cP(Q)$ is the union of the sets of
lattice points of $l\hspace{1pt}Q$ for $0\leq l\leq k$ (for $l=0$ it is just one point, the vertex of the
pyramid). Therefore,
$$\dim(\cC_{k\hspace{1pt}\cP(Q)})=\sum_{l=0}^k\dim(\cC_{l\hspace{1pt}Q}).$$
More on counting lattice points in polytopes see, for example,  \cite{BeRo}.

\section{Application}

We will now apply the results of the previous section to computing the minimum distance
of toric codes whose polytopes are obtained by taking the product, and taking a 
multiple of a pyramid over a polytope.

\begin{Th}\label{T:mix} 
Let $P\subseteq K_q^n$ 
be obtained from $\{0\}$ by successive applications of~$n$ operations of the
following two types
\begin{enumerate}
\item multiplication by a lattice segment  $[0,a_i]$,
\item taking the unit pyramid and scaling the result by $k_i$,
\end{enumerate}
where the first step is {\rm (1)} with $i=1$ and every next step is either {\rm (1)} or {\rm (2)} with the 
{\it next} value of $i$. Let the set of the $a_i$  be indexed by $I\subset \{1,\dots, n\}$ and the set of 
the~$k_i$ be indexed by $J=\{1,\dots, n\}\setminus I$.
Then the minimum distance for the toric code $\cC_P$ equals
\begin{equation}\label{e:mix}
d(\cC_P)=(q-1)^{|J|}\prod_{i\in I}\Big(q-1-a_i\!\!\prod_{j\in J,\, j>i}\! k_j\Big).
\end{equation}

\end{Th}

\begin{pf} The proof is by induction on $n$. If $n=1$ then $I=\{1\}$, $J=\varnothing$ and 
$P$ is a lattice segment of length $a_1$. In this case $\cL(P)$ is the space of univariate polynomials
of degree at most $a_1$, and so $Z_P=a_1$. Therefore, $d(\cC_P)=q-1-a_1$.

Now suppose $Q\subset K_q^{n-1}$ is obtained by applying $n-1$ operations as above.
Then 
$$d(\cC_Q)=(q-1)^{|J'|}\prod_{i\in I'}\Big(q-1-a_i\!\!\prod_{j\in J',\, j>i}\! k_j\Big),\quad |I'|+|J'|=n-1.$$

If at the $n$-th step operation (1) is applied then $P=Q\times [0,a_n]$ and so $I=I'\cup\{a_n\}$, $J=J'$. 
In this case
$$d(\cC_P)=(q-1)^{|J|}\prod_{i\in I'}\Big(q-1-a_i\!\!\prod_{j\in J,\, j>i}\! k_j\Big)(q-1-a_n)$$
by \rt{product}, and \re{mix} follows.

If at the $n$-th step operation (2) is applied then $P=k_n\cP(Q)$ and so $I=I'$, $J=J'\cup\{k_n\}$. 
Note that $P$ is the pyramid of height $k_n$ over the multiple $k_nQ$. Let $m$ be the largest index
in $J$ and suppose
$$\{a_1,\dots, k_m,a_{m+1},\dots,a_{n-1}\}$$ 
is the collection of positive integers used in constructing $Q$. Then we  can obtain $k_nQ$ by applying the same sequence of operations as for $Q$ with a modified
collection 
$$\{a_1,\dots, (k_mk_n),(a_{m+1}k_n),\dots,(a_{n-1}k_n)\}.$$
By induction the minimum 
distance of $\cC_{k_nQ}$ equals
$$d(\cC_{k_nQ})=(q-1)^{|J|-1}\prod_{i\in I,\, i<m}\Big(q-1-a_i\!\!\prod_{j\in J',\, j>i}\! k_jk_n\Big)
\prod_{i=m+1}^{n-1}\left(q-1-a_ik_n\right),$$
which together with \rt{pyramid} gives the required formula \re{mix}. 
\end{pf}

Here are a few immediate corollaries from this theorem. 

\begin{Cor}\label{C:simplex} 
Let $P$ be the $k$-dilate of the standard simplex $P=k\Delta\subseteq K_q^n$.
Then $$d(\cC_P)=(q-1)^{n-1}(q-1-k).$$
\end{Cor}

\begin{Cor}\label{C:box} 
Let $P$ be the rectangular box $P=\prod_{i=1}^n[0,a_i]\subseteq K_q^n$.
Then $$d(\cC_P)=\prod_{i=1}^n(q-1-a_i).$$
\end{Cor}

\begin{Cor}\label{C:simplex-box} 
Let $P$ be the $k$-dilate of the $l$-fold pyramid over the $m$-dimensional unit cube
in $K_q^n$, where $l+m=n$. 
Then $$d(\cC_P)=(q-1)^{l}(q-1-k)^m.$$
\end{Cor}

Using slightly modified arguments we can also compute the minimum distance for the
dilated cross polytope.

\begin{Cor}\label{cross} 
Let $P\subset K_q^n$ be the $k$-dilate of a translated $n$-dimensional cross polytope, 
$P=k\left(\Diamond+(1,\dots,1)\right)$, where 
$$\Diamond=\chull\{e_i,-e_i\ |\ 1\leq i\leq n\}.$$
Then $$d(\cC_P)=(q-1)^{n-1}(q-1-2k).$$
\end{Cor}

\begin{pf} The cross polytope $\Diamond$ is obtained by successively taking the double pyramid, starting with a segment $I$ of length 2. Now the formula follows from the fact that $d(\cC_{kI})=q-1-2k$ and
\rr{doublepyr}, by induction.

\end{pf}

\begin{Rem}\label{R:compare} 
The result of \rc{simplex} follows almost immediately from \cite{Serre}.
Also Corollaries \ref{C:simplex}  and~\ref{C:box} were obtained by  Little and
Schwarz (see \cite{LiSche}) using an analog of Vandermonde matrix. Also
the result of Corollary~\ref{C:box} was obtained by Diego Ruano in~\cite{Ru}.
\end{Rem}

\section{A note on Parameters}\label{S:param}

Recall that the {\it relative minimum distance}  and the {\it information rate} of a code
are defined as  $d(\cC_P)/N$ and $\dim(\cC_P)/N$, where $N$ is the block length of the code.
Furthermore, a family of codes (with $N$ as a parameter) is {\it good} if both 
$d(\cC_P)/N$ and $\dim(\cC_P)/N$ approach positive constants, as $N\to\infty$.
For the toric codes considered in \rt{mix}, we have  $N=(q-1)^n$, 
$\dim(\cC_P)=|P\cap\Z^n|$ (the number of lattice points of $P$), 
and $d(\cC_P)$ computed as in \re{mix}, so it makes sense to consider an infinite
family of toric codes corresponding to $n$-dimensional polytopes $P$ as in \rt{mix}, as $n\to\infty$. 
It turns out that this construction does not produce good families of codes.
We will start with an example.

\begin{Ex} The simplex and the rectangular box do not give a good family of codes.
Indeed, for $P=k\Delta\subseteq K_q^n$ we have a family of toric codes $\cC_P$ with
$N=(q-1)^n$, $\dim(\cC_P)={n+k\choose k}$, and $d(\cC_P)=(q-1)^{n-1}(q-1-k)$,
by \rc{simplex}. Therefore,

$$\frac{d(\cC_P)}{N}=1-\frac{k}{q-1},\quad \frac{\dim(\cC_P)}{N}\sim\frac{n^k}{k!(q-1)^n}\to 0, 
\text{ as }n\to\infty.$$
Similarly, for a rectangular box $P=\prod_{i=1}^n[0,a_i]$, 
we have  $\dim(\cC_P)=\prod_{i=1}^n(a_i+1)$ and 
$d(\cC_P)=\prod_{i=1}^n(q-1-a_i)$, by \rc{box}. Therefore,
$$\frac{d(\cC_P)}{N}=\prod_{i=1}^n\left(1-\frac{a_i}{q-1}\right)\leq
\left(1-\frac{1}{q-1}\right)^n\to 0,\text{ as }n\to\infty.$$
\end{Ex}

Now we consider the  general construction of
\rt{mix}. We fix a subset $I\subseteq\N$ and a 
sequence $\{c_i\ |\ i\in\N\}$ of positive integers $1\leq c_i\leq q-2$ where
$c_i=a_i$ for $i\in I$ and $c_i=k_i$ for $i\in\N\setminus I$.
Assume that for every $n\in\N$ the $n$-dimensional polytope $P_n$ constructed as in \rt{mix} 
using the truncated sequence $\{c_i\ |\ 1\leq i\leq n\}$ lies
in the $n$-cube $K_q^n$. Then we obtain a family of toric codes $\cC_n$ associated to
the family of polytopes  $P_n$. The next proposition shows that no such family is a good family of codes.

\begin{Prop}  Let $C_n$ be a family of toric codes constructed as above.
Then either $d(\cC_n)/N$ or $\dim(\cC_n)/N$ approaches zero, as $n\to\infty$.
\end{Prop}

\begin{pf} Let $\{c_i\ |\ i\in\N\}$ be a sequence producing a family
of polytopes $P_n$, as above.
First note that since $P_n$ lies in $K_q^n$, its linear dimensions are bounded. 
This implies that $k_i=1$ for all sufficiently large values of $i$. Furthermore, we will
show that for $d(\cC_n)/N$ to approach a positive constant, the subset $I\subseteq\N$ must
be finite. Indeed, for every polytope $P_n$ in the family, by \re{mix}, we have 
$$\frac{d(\cC_n)}{N}=\prod_{i\in I_n}\Big(1-\frac{a_i}{q-1}\!\!\prod_{j\in J_n,\, j>i}\! k_j\Big),$$
where $I_n=I\cap\{1,\dots,n\}$ and $J_n=\{1,\dots,n\}\setminus I_n$. Therefore, using 
 $a_i\geq 1$ for all $i$, we obtain
\begin{equation}\label{e:bounds}
\frac{d(\cC_n)}{N}\leq\prod_{i\in I_n}\Big(1-\frac{a_i}{q-1}\Big)\leq
\left(1-\frac{1}{q-1}\right)^{|I_n|}.
\end{equation}
If $I\subseteq\N$ is infinite then
$|I_n|\to\infty$, as $n\to\infty$, which implies that the right hand side of \re{bounds}, and 
hence so does $d(\cC_n)/N$, approaches zero, as $n\to\infty$.

Now if $I$ is finite and the $k_i$ stabilize to 1 then for all large enough $n$
the polytope $P_n$ is obtained simply by building  the unit pyramid over the previously
constructed polytope $P_{n-1}$. Therefore $\dim(\cC_n)=\dim(\cC_{n-1})+1$, for large enough $n$.
This implies that $\dim(\cC_n)/N\to 0$, as $n\to\infty$.
\end{pf}

It would be interesting to find an infinite good family of toric codes as $n\to\infty$.


\section{Examples}\label{S:ex}

In this section we look at a few examples of  toric codes of low dimension and 
compare their parameters. We will start with a 2-dimensional example.

\begin{Ex}\label{ex:1} 
Let $Q$ be the lattice triangle with vertices $(1,0)$, $(0,3)$, and $(3,1)$ (see
\rf{triangle}), and let $\cC_Q$ be the corresponding toric code.
 \begin{figure}[h]
\centerline{
 \scalebox{0.6}
 {
\input{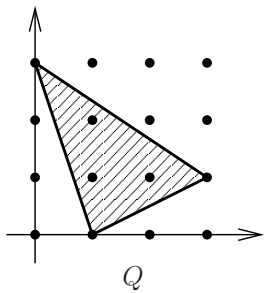}}}
\caption{Tirangle with 6 lattice points.}
\label{F:triangle}
\end{figure}
 The triangle has~6 lattice points, so $\dim(\cC_Q)=6$.
Note that the polynomial $xy(x-a)(x-b)$, for $a,b\in\F_q$, $a\neq b$, lies in $\cL(Q)$
and has weight $(q-1)(q-3)$, so the minimum distance satisfies $d(\cC_Q)\leq(q-1)(q-3)$.
In fact, according to Theorem 2.6 of \cite{SoSo}, $d(\cC_Q)=(q-1)(q-3)$ for all $q\geq 37$. 
Combining this with Magma \cite{magma} computations for small $q$ we obtain the following answer:
$$d(\cC_Q)=\begin{cases}(q-1)(q-3), & \text{ if }q\neq 8\\ \,28, & \text{ if }q=8.\end{cases}$$
Note that for $q=5$ we obtain a $[16, 6, 8]$-code, which is a best known code over~$\F_5$ (see \cite{table}).
\begin{Rem}
Here is a reason why $q=8$ is a special case. Notice 
that $\cL(Q)$ contains a polynomial $f(x,y)=x^3y+y^3+x$. Let $Y$ be the closure 
of $\{f=0\}$ in the toric surface defined by~$Q$. Then $Y$ produces an example
of a genus 3 absolutely irreducible curve over
$\F_8$ on which the upper part of Serre's Bound is attained: The 
the number of $\F_8$-rational points on $Y$
is $$q+1+g\lfloor{2\sqrt{q}}\rfloor=8+1+3\cdot\lfloor{2\sqrt{8}}\rfloor=24$$ (see \cite{Sti} Example VI.3.8). Three of the rational points on $Y$ lie
``at infinity'' (on the invariant orbits of the toric surface), hence $Y$
has $21$ rational points in the torus $(\F_8^*)^2$. This implies that the weight of $f$ is
$w(f)=49-21=28$ and it is the minimum weight here.
\end{Rem}
 \end{Ex}
 
 Now we will consider examples of 3-dimensional toric codes related to the
 toric surface code in \rex{1}.
 
\begin{Ex} \label{ex:2}
Consider the prism $R=Q\times[0,1]$ with base $Q$ from \rex{1}. The toric code $\cC_R$ has
parameters $N=(q-1)^3$, $\dim(\cC_R)=12$, and, for $q\neq 8$,
$$d(\cC_R)=d(\cC_Q)(q-2)=(q-1)(q-2)(q-3),$$
by \rt{product}. In particular, over $\F_5$ we obtain a $[64, 12, 24]$-code. According to the table of
linear codes \cite{table}, the best known linear code over $\F_5$ with $N=64$ and $k=12$ has  
$d=36$, which by far exceeds the minimum distance of $\cC_R$.
\end{Ex}

\begin{Ex} \label{ex:3}
Let $\cP(Q)$ be the unit pyramid over $Q$ from \rex{1}. The toric code $\cC_{P(Q)}$ has
parameters $N=(q-1)^3$, $\dim(\cC_{P(Q)})=7$, and, for $q\neq 8$, 
$$d(\cC_{P(Q)})=d(\cC_Q)(q-1)=(q-1)^2(q-3),$$
by \rt{pyramid}. In particular, over $\F_5$ we obtain a $[64, 7, 32]$-code. Again, for a code with
$N=64$ and $d=32$, its dimension is very small. The best known code with $N=64$ 
and $d\geq 32$ has dimension $k=14$ (see \cite{table}).
\end{Ex}

The above two examples suggest that one should look for codes whose polytopes are not
products of or pyramids over smaller dimensional polytopes. Below we show one
example of a 3-dimensional toric code with parameters $[64,13, 31]$, which is much closer
to the best known $[64,14, 33]$-code than the ones in Examples \ref{ex:2} and \ref{ex:3}.

\begin{Ex} \label{ex:4} Let $P$ be the convex hull of five lattice points: $(0,3,0)$, $(1,0,0)$, $(3,1,0)$,
$(1,1,2)$, and $(2,3,3)$. The polytope $P$ contains $Q$ as its base (see \rf{triangle}),
so $d(\cC_P)\leq d(\cC_Q)(q-1)$. For $q=5$ Magma produces $d(\cC_P)=31$.
Furthermore, $P$ has ~13 lattice points: five vertices, three more in the base, and 
$(1,1,1)$, $(2,1,1)$, $(1,2,1)$, $(2,2,1)$ and $(2,2,2)$; hence,
$\dim(\cC_{P})=13$. Therefore $\cC_P$ is a  $[64, 13, 31]$-code over $\F_5$.
\end{Ex}



\begin{thebibliography}{99}


\bibitem[1]{BeRo} M. Beck, S. Robins, {\em Computing the continuous discretely. Integer-point enumeration in polyhedra.} Undergraduate Texts in Mathematics. Springer, New York, 2007.

\bibitem[2]{magma}W.~Bosma, J.~Cannon, and C.~Playoust,
{\em The Magma algebra system. I. The user language.} 
J. Symbolic Comput., 24(3-4):235-265, 1997 

\bibitem[3]{F} W. Fulton, {\em Introduction to Toric Varieties},
Princeton Univ. Press, Princeton, 1993


\bibitem[4]{table} Markus Grassl,
{\em Bounds on the minimum distance of linear codes and quantum codes},
Online available at {\tt http://www.codetables.de},
Accessed on 2010-04-06.

\bibitem[5]{Ha1} J. Hansen, {\em Toric Surfaces and Error--correcting
    Codes} in Coding Theory, Cryptography, and Related Areas, Springer
  (2000), pp. 132-142.
  
\bibitem[6]{Ha2} J. Hansen, Toric varieties Hirzebruch surfaces and error-correcting codes, Appl. Algebra Engrg. Comm. Comput. 13 (2002), pp. 289Ð300
  
\bibitem[7]{Jo} D. Joyner, {\em Toric codes over finite fields}, Appl.
Algebra Engrg. Comm. Comput.,
15 (2004), pp. 63--79.
  
\bibitem[8]{LiSche} J. Little, H. Schenck,
{\em Toric Surface Codes and Minkowski sums},
SIAM J. Discrete Math.  20  (2006),  no. 4, 999--1014.

\bibitem[9]{LiSchw} J. Little, R. Schwarz, 
{\em On toric codes and multivariate Vandermonde matrices}, 
Appl. Algebra Engrg. Comm. Comput. 18 (4) (2007), pp. 349--367.

\bibitem[10]{Ru}Diego Ruano, 
{\em On the parameters of $r$-dimensional toric codes},
Finite Fields and Their Applications {\bf 13} (2007), pp. 962--976.

\bibitem[11]{Serre} J.-P. Serre, {\em Lettre a M. Tsfasman}, Ast\'{e}risque {\bf 198-199-200} (1991), pp. 351--353. 


\bibitem[12]{SoSo} I.~Soprunov, J.~Soprunova, {\em Toric surface codes and Minkowski length of polygons}, SIAM J. Discrete Math. {\bf 23}, Issue 1, (2009) pp. 384-400 


\bibitem[13]{Sti} H. Stichtenoth,
{\em Algebraic Function Fields and Codes}, Springer-Verlag, (1993).

\bibitem[14]{TV} M. Tsfasman, S. Vl\u adu\c t, 
{\em Algebraic-geometric codes}, Kluwer, Dordrecht, (1991).
\end{thebibliography}
\end{document}